\documentclass[12pt]{amsart}
\usepackage{amscd,amssymb,verbatim}

\newcommand{\sA}{\mathcal{A}}
\newcommand{\sB}{\mathcal{B}}

\newcommand{\sM}{\mathcal{M}}
\newcommand{\sN}{\mathcal{N}}
\newcommand{\sO}{\mathcal{O}}
\newcommand{\sR}{\mathcal{R}}

\newcommand{\1}{\mathbf{1}}

\newcommand{\C}{\mathbf{C}}
\renewcommand{\P}{\mathbf{P}}
\newcommand{\Q}{\mathbf{Q}}

\newcommand{\bS}{\mathbb{S}}
\newcommand{\Z}{\mathbf{Z}}

\renewcommand{\epsilon}{\varepsilon}

\newcommand{\cf}{{\it cf. }}
\newcommand{\resp}{{\it resp. }}

\newcommand{\End}{\operatorname{End}}
\newcommand{\Hom}{\operatorname{Hom}}

\newcommand{\Aut}{\operatorname{Aut}}

\newcommand{\Trd}{\operatorname{Trd}}
\newcommand{\Nrd}{\operatorname{Nrd}}
\newcommand{\Tr}{\operatorname{Tr}}
\newcommand{\tr}{\operatorname{tr}}
\newcommand{\gr}{{\operatorname{gr}}}
\newcommand{\ord}{{\operatorname{ord}}}
\renewcommand{\int}{{\operatorname{int}}}

\renewcommand{\Vec}{\operatorname{Vec}}
\newcommand{\Rep}{\operatorname{Rep}}

\newcommand{\num}{{\operatorname{num}}}

\newcommand{\iso}{\overset{\sim}{\longrightarrow}}

\swapnumbers
\newtheorem{lemma}{Lemma}[section]
\newtheorem{thm}[lemma]{Theorem}
\newtheorem{cor}[lemma]{Corollary}
\newtheorem{prop}[lemma]{Proposition}
\theoremstyle{definition}
\newtheorem{rk}[lemma]{Remark}

\newtheorem{defn}[lemma]{Definition}
\newtheorem{nota}[lemma]{Notation}

\newcounter{spec}
\newenvironment{thlist}{\begin{list}{\rm{(\roman{spec})}}%
{\usecounter{spec}\labelwidth=20pt\itemindent=0pt\labelsep=10pt}}%
{\end{list}}

\numberwithin{equation}{section}

\begin{document}

\title{On the multiplicities of a motive}
\author{Bruno Kahn}
 \address{Institut de Math\'ematiques de
 Jussieu\\175--179 rue du
 Chevaleret\\ \break 75013
 Paris\\France.}
\email{kahn@math.jussieu.fr}
\date{January 12, 2007}
\begin{abstract} We introduce the notion of multiplicities for an object $M$ in a semi-simple
rigid tensor category $\sA$, as a collection of central scalars which relate the categorical
trace with the ring-theoretic trace. Multiplicities turn out to be rational integers in
important cases, most notably when $\sA$ is of ``homological origin". We show
that this integrality condition has simple consequences, like the rationality and a functional
equation for the zeta function of an (invertible) endomorphism. An example is the category of
pure motives modulo numerical equivalence with rational coefficients over a field $k$; if $k$
is finite and $M$ is of abelian type, its multiplicities are all equal to $\pm 1$.
\end{abstract}
\maketitle

\tableofcontents

\section*{Introduction} The aim of this article, in the spirit of
\cite{nrsm}, is to study abstractly the properties of categories of pure
motives and to make clear(er) which of them are formal and which are of a
more arithmetic-geometric nature. 

We work with a rigid additive tensor category $\sA$ such that $K=\End(\1)$ is a
field of characteristic $0$. We shall be interested in the
\emph{multiplicities} of an object $M\in \sA$: when $\sA$
is semi-simple, they are a collection of central scalars which relates the
categorical trace with the ring-theoretic trace (Proposition \ref{p1}). It
turns out that the condition for these multiplicities to be \emph{integers}
or, better, to be so after extending scalars from $K$ to its algebraic
closure, is very well-behaved and is satisfied in many important cases.
Namely: 

\begin{itemize} 
\item The full subcategory $\sA_\int$ of $\sA$ formed by
such objects is thick, tensor, rigid, contains the ``Schur-finite"
objects (those which are killed by a nonzero Schur functor), and is preserved under tensor
functors to another semi-simple rigid category (Corollary \ref{c3}). 
\item $\sA_\int=\sA$ if $\sA$ is of ``homological origin" (Theorem
\ref{t2}). The category of pure motives over a field modulo numerical
equivalence is semi-simple thanks to Jannsen's theorem \cite{jannsen}, and
of homological origin.  
\end{itemize}

When the multiplicities are integers, we prove that the zeta function of an endomorphism $f$ of
$M$ is rational (with an explicit formula) and satisfies a functional equation if $f$ is
invertible (Theorem \ref{t3}): in the case of motives over a finite field, this shows that
these depend on less than the existence of a Weil cohomology theory. We also get some
elementary cases where homological equivalence equals numerical
equivalence for formal reasons in Proposition \ref{p5} c): of course, this
remains far from leading to a proof of this famous standard conjecture!

In Section \ref{5}, we formulate a version of the Tate
conjecture for motives over a finite field in an abstract set-up. My
initial motivation was to see what multiplicities had to say on this
conjecture; this turns out to be disappointing (see Theorem \ref{t4}
and Remark \ref{r1}) but I find it amusing and perhaps enlightening
that most of its known equivalent versions
carry out in this abstract context: see Theorem \ref{t4} and Corollary
\ref{c1}. The proof of
\cite{cell} that under the Tate conjecture and Kimura-O'Sullivan
finite dimensionality, rational and numerical equivalences agree over
a finite field also carries out abstractly: see Theorem \ref{tcell}.

\subsection*{Acknowledgements} I thank Yves Andr\'e, Alessio del Padrone and Chuck Weibel for helpful exchanges
during the preparation of this paper. Part of this research was done during a visit to the Tata
Institute of Fundamental Research, partially funded by the Indo-French Institute for
Mathematics: I would like to thank both institutions for their hospitality and support.

\enlargethispage*{20pt}

\subsection*{Terminology and notation} Let $\sA$ be a rigid $K$-linear tensor category, where
$K$ is a field of characteristic $0$;  we also assume that $\End(\1)=K$. In the sequel of this
article, we shall abbreviate this by saying that $\sA$ is a
\emph{rigid $K$-category}. Since we shall refer to Deligne's article
\cite{deligne} several times, it is worth stressing that we do not assume $\sA$
abelian, unlike in {\it loc. cit}. We write $\sA^\natural$ for the pseudo-abelian hull
(idempotent completion) of $\sA$.

If $M\in \sA$, we shall say (as has become common practice) that $M$ is \emph{Schur-finite}
if  there exists a nonzero Schur functor $\bS$ such that $\bS(M)=0$ and
\emph{finite-dimensional} (in the sense of Kimura-O'Sullivan) if $M\simeq M_+\oplus M_-$ where
$M_+$ (\resp $M_-$) is killed by some nonzero exterior (\resp symmetric) power functor. We say
that $M_+$ is \emph{positive} and $M_-$ is \emph{negative}. It is known that finite-dimensional
implies Schur-finite (\cf \cite[1.7]{deligne}). For properties of finite-dimensional objects
(\resp of Schur functors) we refer to \cite{kim} and \cite[\S 9]{nrsm} (\resp to
\cite{deligne}).

\section{Multiplicities in semi-simple rigid tensor categories}\label{s1}

Let $M\in \sA$. The \emph{trace} of an endomorphism $f\in \End(M)$ is the
element $\tr(f)\in \End(\1)=K$ defined by the composition 
\[\begin{CD}
\1@>\eta>> M^*\otimes M@>1\otimes f>> M^*\otimes M@>R>> M\otimes
M^*@>\epsilon>> \1 
\end{CD}\] 
where $R$ is the switch and $\eta,\epsilon$
are the duality structures of $M$. 

\subsection*{Special case} We shall denote the trace of $1_M$ by
$\chi(M)$ and call it the \emph{Euler characteristic of $M$}.

The trace is $K$-linear and has the
following properties:  
\begin{equation}\label{eq4} \tr(fg) = \tr(gf),
\quad \tr(f\otimes g) = \tr(f)\tr(g), \quad \tr({}^tf)=\tr(f). 
\end{equation}

Suppose that $\sA$ is semi-simple. Then $\End_\sA(M)$ is a semi-simple
$K$-algebra, hence has its own trace, and we want to compare the
categorical trace with the ring-theoretic trace. We normalise conventions
as follows: 

\begin{defn}\label{d1} a) Let $A$ be a finite-dimensional simple
$K$-algebra. We write:  
\begin{itemize} 
\item $Z(A)$ for the centre of $A$;  
\item $\delta(A)=[Z(A):K]$;  
\item $d(A)=[A:Z(A)]^{1/2}$. 
\end{itemize} 
We define the \emph{reduced trace} of $A$ as \[\Trd_A =
\Tr_{Z(A)/K}\circ \Trd_{A/Z(A)}.\] If $A=\prod A_i$ is semi-simple, with
simple components $A_i$, we define $\Trd_A := \sum_i \Trd_{A_i}$.\\ b) If
$A=\End_\sA(M)$, we set 
\begin{itemize} 
\item $Z_i(M)=Z(A_i)$;  
\item $\delta_i(M)=\delta(A_i)$;  
\item $d_i(M)=d(A_i)$;  
\item $\Trd_M=\Trd_A$. 
\end{itemize} 
\end{defn}

\begin{prop} \label{p1} There exists a unique element $\mu(M)\in \End(M)$
such that
\[\tr(f)= \Trd_M(\mu(M)f)\]
for any $f\in \End(M)$. Moreover, $\mu(M)$ is central and invertible.
Hence, if $(e_i)$ denotes the set of central idempotents of $A=\End(M)$
corresponding to its simple factors $A_i$, we may write
\[\mu(M)=\sum_i \mu_i(M) e_i\]
with $\mu_i(M)\in Z_i(M)$.
\end{prop}

\begin{proof} Since $\End(M)$ is semi-simple, $(f,g)\mapsto \Trd_M(fg)$
is nondegenerate,
which proves the existence and uniqueness of $\mu(M)$. Moreover,
\[\Trd_M(\mu(M)fg)=\tr(fg)=\tr(gf)=\Trd_M(\mu(M)gf)=\Trd_M(f\mu(M)g)\]
and the non-degeneracy also yields the centrality of $\mu(M)$. This
element is invertible because the ideal $\sN$ is $0$ for $\sA$
\cite[7.1.7]{nrsm}. The last assertion is obvious.
\end{proof}

\begin{lemma}\label{l9} a) We have $\mu(M^*) = {}^t\mu(M)$.\\
b) Suppose that $K$ is algebraically closed and $M$ is simple. Then
$\mu(M) = \chi(M)$.
\end{lemma}

\begin{proof} a) follows easily from \eqref{eq4} and the fact that the
transposition
induces an anti-isomorphism from $\End(M)$ onto $\End(M^*)$. b) is
obvious, since then $\End(M)
= K$ (recall that, by definition, $\chi(M)=\tr(1_M)$).
\end{proof}

\begin{rk}\label{l1} If $\sA$ is pseudo-abelian (hence abelian, \cite[Lemma 2]{jannsen}), the
idempotents $e_i$ of
Proposition \ref{p1} yield the decomposition $M=\bigoplus M_i$ of $M$
into its \emph{isotypical
components}. In particular, if $S$ is simple, then $\mu(S^n)=\mu(S)$ for
any $n\ge 1$.\\ 
On the other hand, it is difficult to relate $\mu(M_1),\mu(M_2)$ and
$\mu(M_1\otimes M_2)$ in general
because it is difficult to say something of the map $\End(M_1)\otimes_K
\End(M_2)\to\End(M_1\otimes M_2)$: it is not even true in general that
such a homomorphism sends the centre into the centre. For the same reason,
it is difficult to state general facts on the behaviour of the invariant
$\mu$ under tensor functors. We shall see that the
situation improves considerably in the case of \emph{geometrically
integral type}, discussed in the next section.
\end{rk} 

\section{Integral multiplicities}

In all this section, $\sA$ is a semi-simple rigid $K$-category.

\begin{defn} a) An object $M\in\sA$ is \emph{of integral type} if 
the scalars $\mu_i(M)$ of Proposition \ref{p1} belong to $\Z$.\\
b) $M$ is \emph{geometrically of
integral type} if $M_{\bar K}\in \sA_{\bar K}$ is of integral type, where
$\bar K$ is an
algebraic closure of $K$.\\
c) $\sA$ is \emph{of integral type} (\resp \emph{geometrically of
integral type}) if every
$M\in \sA$ is of integral type (\resp geometrically of integral type).  
\end{defn}

\begin{prop}\label{p2} a) If $M$ is of integral type, we have 
\begin{equation}\label{eq3}
\mu_i(M) =\frac{\tr(e_i)}{\delta_i(M)d_i(M)}
\end{equation} 
for any $i$.\\ 
b) Direct sums and direct summands of objects of integral type are of
integral type. Similarly
for  geometrically of integral type. In particular, $\sA$ is of integral
type (\resp geometrically of integral type) if and only if its
pseudo-abelian envelope is.\\   
c) If $M$ is
geometrically of integral type, then it is of integral type. Moreover, if
this is the case, the
invariants
$\mu_i(M)$ are ``geometric" in the sense that if $L/K$ is any extension,
then $\mu_i(M) = \mu_{i,\alpha}(M_L)$ for any simple factor $A_{i,\alpha}$
of $A_i\otimes_K L$.\\
d) $M\in \sA$ is geometrically of integral type if and only if, in
$\sA_{\bar K}^\natural$,
the Euler characteristic of every simple summand of $M_{\bar K}$ is an integer.\\
e)  If $M$ is Schur-finite, it is
geometrically of integral type.\\ 
f) If $M$ is (geometrically) of integral type, so is $M^*$.
\end{prop}

\begin{proof} a) and b) are obvious. For c), we have
the decomposition
\[Z_i(M)\otimes_K \bar K\iso \prod_\alpha \bar K\]
where $\alpha$ runs through the distinct $K$-embeddings of $Z_i(M)$ into
$\bar K$.
Correspondingly, $A_i\otimes_K {\bar K}$ decomposes as a direct product
\[A_i\otimes_K {\bar K} \simeq \prod_\alpha A_i^\alpha\]
with $A_i^\alpha$ simple over $\bar K$. This gives a decomposition
\[e_i\otimes_K 1=\sum_\alpha e_i^\alpha\]
into central idempotents. But clearly, $\mu(M_{\bar K})=\mu(M)\otimes_K
1$. By hypothesis,
the images of $\mu_i(M)$ in $\bar K$ under the embeddings $\alpha$ are
rational integers, which
implies that $\mu_i(M)$ is itself a rational integer. The additional
claim of c) immediately follows from this proof.

d) follows immediately from Lemma \ref{l9} a). 

For e), if $M$ is Schur-finite, so is $M_{\bar K}\in \sA_{\bar K}^\natural$; all simple direct
summands of $M_{\bar K}$ are Schur-finite as well, hence their Euler characteristics are
rational integers. This immediately follows from the main result of  \cite{deligne}, but one
can more elementarily use Proposition 2.2.2 of A. del Padrone's thesis \cite{delp}, which
generalises the case of finite-dimensional objects \cite[7.2.4 and 9.1.7]{nrsm}. The conclusion
now follows from d).

Finally, f) follows from Lemma \ref{l9} b).
\end{proof}

\begin{rk} C. Weibel raised the question whether the converse of e) is true. I don't know any
counterexample; it holds at least if $\sA_\int$ is of homological origin in the sense of
Definition \ref{d3} b) (see Theorem \ref{t2}).
\end{rk}

\begin{thm}\label{T1} Let
$M,N\in\sA$ be geometrically of integral type, $(e_i)$ the central
idempotents of $\End(M)$ and
$(f_j)$ the central idempotents of $\End(N)$. For a pair $(i,j)$, let
$A_{ij}$ be the
semi-simple algebra $(e_i\otimes f_j) \End(M\otimes N)(e_i\otimes
f_j)$. Then one has formulas of the type
\[\mu_i(M)\mu_j(N) = \sum_k m_k \mu_k(M\otimes N)\]
where $k$ indexes the simple factors  of $A_{ij}$ and the $m_k$ are
integers $\ge 0$.
Moreover, for any $k$, there is such a formula with $m_k>0$.\\
In particular, $M\otimes N$ is geometrically of integral type.
\end{thm}

\begin{proof} We proceed in 2 steps:

1) $\End(M)$ and $\End(N)$ are split. By Proposition \ref{p2} b), we may
assume that
$\sA$ is pseudo-abelian. This allows us to assume $M$ and $N$
\emph{simple}, hence
$\End(M)=\End(N)=K$ and $A_{ij}=\End(M\otimes N)$. 
Using Formula \eqref{eq3} to compute $\tr(1_{M}\otimes 1_{N})$ in two
different ways, we get
the formula
\begin{equation}\label{eq2}
\mu(M)\mu(N)= \sum m_k\mu_k(M\otimes N)
\end{equation}
with $m_k = \delta_k(M\otimes N)d_k(M\otimes N)$.

Coming back to the case where $\sA$ is not necessarily
pseudo-abelian and $M,N$ not necessarily simple, this gives the
formula
\[m_k = \delta_k(A_{ij})\frac{d_k(A_{ij})}{d_i(M)d_j(N)}\]
(see Remark \ref{l1}), and the previous argument shows ungrievously that this is an integer.

2) The general case. Extending scalars to $\bar K$ and using Proposition
\ref{p2} c), we are
reduced to 1) as follows: for any $\alpha:Z_i(M)\to \bar K$ and any
$\beta:Z_j(M)\to \bar
K$, we have a formula with obvious notation:
\[\mu_i^\alpha(M_{\bar K})\mu_j^\beta(N_{\bar K})= \sum_k\sum_\gamma
m_k^\gamma\mu_k^\gamma((M\otimes N)_{\bar K})
\]
where, for each $k$, $\gamma$ runs through the embeddings of
$Z_k(M\otimes N)$ into $\bar
K$. By Remark \ref{l1}, this gives a formula as wanted.

It remains to prove that, given a simple
factor $A_k$ of $A_{ij}$, one may find a formula with $m_k>0$. For this, it suffices
to show that there is a pair $(\alpha,\beta)$ such that 
\[\Hom_{\bar K}(A_k\otimes_K \bar K,(e_i^\alpha\otimes f_j^\beta)
(A_{ij}\otimes_K \bar
K)(e_i^\alpha\otimes f_j^\beta))\ne 0.\]

This is obvious, since $\Hom_K(A_k,A_{ij})\ne 0$ and $A_{ij}\otimes_K \bar
K=\prod_{\alpha,\beta} (e_i^\alpha\otimes f_j^\beta) (A_{ij}\otimes_K
\bar K)(e_i^\alpha\otimes
f_j^\beta))$.
\end{proof}

\begin{cor}\label{c2} Assume that $M$ and $N$ are simple and that, in
Theorem \ref{T1}, all
terms
$\mu_k(M\otimes N)$ have the same sign. Then we have $|\mu_k(M\otimes
N)|\le |\mu(M)\mu(N)|$ for
all
$k$. If $|\mu(M)|=|\mu(N)|=1$, then $A=\End(M\otimes N)$ is ``geometrically
simple" in the sense that $A\otimes_K \bar K$ is a matrix algebra over
$Z(M\otimes
N)\otimes_K \bar K$ (otherwise said, $A$ is an Azumaya algebra over its
centre). Moreover,
$\mu(M\otimes N)=\mu(M)\mu(N)$.
\end{cor}

\begin{proof} This follows from the last statement of Theorem \ref{T1}. In the
case where $|\mu(M)|=|\mu(N)|=1$, Formula \eqref{eq2} gives the conclusion.
\end{proof}

\begin{cor}\label{c3} a) The full subcategory $\sA_\int$ of $\sA$
consisting of geometrically
integral objects is a thick rigid tensor subcategory of $\sA$ containing the Schur-finite objects.\\
b) Let $F:\sA\to \sB$ be a $\otimes$-functor to another rigid semi-simple
$K$-category. Then
$F(\sA_\int)\subseteq \sB_\int$.
\end{cor}

\begin{proof} a) follows from Proposition \ref{p2} and Theorem \ref{T1}. b) follows from
Proposition \ref{p2} d).
\end{proof}

\section{Application: the zeta function of an endomorphism}

\begin{defn} Let $\sA$ be a rigid $K$-category, $M\in \sA$ and $f\in
\End(M)$. The \emph{zeta
function} of $f$ is
\[Z(f,t) = \exp \left(\sum_{k\ge 1} \tr(f^n) \frac{t^n}{n}\right)\in K[[t]].\]
\end{defn}

\begin{thm}\label{t3} Suppose that $\sA$ is semi-simple and that $M\in
\sA$ is of integral type.
Then,\\ 
a) For any $f\in \End(M)$, $Z(f,t)\in K(t)$. More precisely, one has with
the notation
of Definition \ref{d1}
\[Z(f,t) = \prod_i \Nrd_{A_i}(e_i-e_i f t)^{-\mu_i(M)}\]
where, for all $i$, $\Nrd_{A_i}(e_i-e_i f t)  
:=N_{Z_i(M)/F}\Nrd_{A_i/Z_i(M)}(e_i-e_i ft)$ denotes the inverse reduced
characteristic
polynomial of the element $e_i f$ if $A_i$.\\
b) If $f$ is invertible, one has the functional equation
\[Z(f^{-1},t^{-1})= (-t)^{\chi(M)} \det(f) Z(f,t)\]
where $\chi(M)=\tr(1_M)$ and $\det(f) =\prod_i\Nrd_{A_i}(e_if)^{\mu_i(M)}$.
\end{thm}

%\enlargethispage*{20pt}

\begin{proof} a) Applying the formula of Proposition \ref{p1}, we get
\begin{multline*}
Z(f,t) = \exp \left(\sum_{k\ge 1} \Trd_M(\mu(M)f^n) \frac{t^n}{n}\right)\\
= \exp \left(\sum_{k\ge 1}\sum_i \Trd_M(\mu_i(M)e_if^n) \frac{t^n}{n}\right)\\
= \prod_i\exp \left(\sum_{k\ge 1} \Trd_{A_i}((e_if)^n)
\frac{t^n}{n}\right)^{\mu_i(M)}\end{multline*} 
and the conclusion follows from the well-known
linear algebra identity
\[\exp \left(\sum_{k\ge 1} \Trd_{A_i}((e_if)^n)
\frac{t^n}{n}\right)=\Nrd_{A_i}(e_i-e_i f t)^{-1}.\]

For b), we write
\[\Nrd_{A_i}(e_i-e_if^{-1}t^{-1})= 
\Nrd_{A_i}(-e_if^{-1}t^{-1})\Nrd_{A_1}(e_i-e_ift)\]
hence
\begin{multline*}
Z(f^{-1},t^{-1})=\prod_i\Nrd_{A_i}(e_i-e_if^{-1}t^{-1})^{-\mu_i(M)}\\
=\prod_i
\Nrd_{A_i}(-e_if^{-1}t^{-1})^{-\mu_i(M)}\Nrd_{A_1}(e_i-e_ift)^{-\mu_i(M)}\\
= \prod_i\Nrd_{A_i}(-e_if^{-1}t^{-1})^{-\mu_i(M)} Z(f,t)
\end{multline*}
and
\begin{multline*}
\prod_i\Nrd_{A_i}(-e_if^{-1}t^{-1})^{-\mu_i(M)}=\\
(-t)^{\sum_i\mu_i(M)d_i(M)\delta_i(M)}\prod_i\Nrd_{A_i}(e_if)^{\mu_i(M)}=
(-t)^{\chi(M)}\det(f).
\end{multline*}
\end{proof}

\begin{rk} The definition of $\det$ shows that
\[\det(1-ft)= Z(f,t)^{-1}\]
if the left hand side is computed in $\sA_{K(t)}$.
\end{rk}

\section{Multiplicities in rigid tensor categories of homological type}

\begin{defn}\label{d3} a) A rigid $K$-category $\sA$ is \emph{of
homological type} if there exists a tensor functor
\[H:\sA\to \Vec^\pm_L\]
where $L$ is an extension of $K$ and $\Vec^\pm_L$ is the tensor category
of $\Z/2$-graded
finite-dimensional $L$-vector spaces, provided with the Koszul rule for
the commutativity
constraint. We say that $H$ is a \emph{realisation} of $\sA$\footnote{When $\sA$ is abelian, this is what Deligne calls a super-fibre functor in \protect\cite{deligne}, except that we do not require any exactness or faithfulness property here.}. \\ 
We say that $\sA$ is \emph{neutrally of homological type} if one may
choose $L=K$.
b) A semi-simple rigid $K$-category $\bar \sA$ is \emph{of homological
origin} (\resp
\emph{neutrally of homological origin}) if it is $\otimes$-equivalent to
$\sA/\sN$, where $\sA$
is a rigid $K$-category of homological type (\resp neutrally of
homological type) and
$\sN=\sN(\sA)$ is the ideal of morphisms universally of trace $0$.
\end{defn}

\begin{lemma}\label{l6}  If $\sA$ is of homological type, $\sA/\sN$ is
semi-simple. If moreover
it is neutrally of homological type and the corresponding realisation $H$
is faithful, the
functor
$\sA\to \sA/\sN$ has the idempotent lifting property.
\end{lemma}

\begin{proof} The first statement follows from \cite[Th. 1 a)]{err-nrsm}.
For the second, let
$M\in \sA$ and $\bar M$ its image in $\bar \sA$. The hypothesis implies
that $\End_\sA(M)$ is a
finite-dimensional $K$-algebra. Let $\sR$ be its radical: it is nilpotent
and
contained in $\sN(M,M)$ by \cite[Th. 1 a)]{err-nrsm}. Thus $\End_{\bar
\sA}(\bar M)$ is a
quotient of the  semi-simple algebra $\End_\sA(M)/\sR$. Therefore we may
lift orthogonal
idempotents of $\End_{\bar \sA}(\bar M)$ to orthogonal idempotents of
$\End_\sA(M)$, first in
$\End_\sA(M)/\sR$ and then in $\End_\sA(M)$ itself.
\end{proof}

\begin{nota} Let $\sA$ be of homological type. For $M\in \sA$, we
  write $\delta_i(M),d_i(M),\mu_i(M)$ for $\delta_i(\bar M),d_i(\bar
  M),\mu_i(\bar M)$, where $\bar M$ is the image of $M$ in $\bar\sA$.
\end{nota}

\begin{lemma}\label{l5} Let $E$ be an extension of $K$. If $\bar\sA$ is
of homological origin,
then
$\bar \sA_E:=\bar\sA\otimes_K E$ is also of homological origin.
\end{lemma}

\begin{proof} Let $\sA$ of homological type be such that $\sA/\sN\simeq
\bar \sA$, and let
$H:\sA\to \Vec_L^\pm$ be a realisation of $\sA$. Consider the tensor functor
\[H_E:\sA_E\to \Vec_{L\otimes_K E}^\pm\]
given by $H_E(M) = H(M)\otimes_K E$.
Here $L\otimes_K E$ is not a field in general, but we can map it to one
of its residue fields
$L'$. Then the composite functor
\[H':\sA_E\to \Vec_{L'}^\pm\]
is a tensor functor. To conclude, it suffices to observe that $\bar
\sA_E\simeq
\sA_E/\sN(\sA_E)$ by \cite[Lemme 1]{err-nrsm}.
\end{proof}

\begin{lemma}\label{l7} Suppose that $\bar \sA$ is neutrally of
homological origin. Then the
pseudo-abelian envelope of $\bar\sA$ is also neutrally of homological origin.
\end{lemma}

\begin{proof} A realisation $H$ with coefficients $K$ extends to the
pseudo-abel\-ian envelope
$\sA^\natural$ of $\sA$, since $\Vec_K^\pm$ is pseudo-abelian. On the
other hand, Lemma
\ref{l6} implies that $\sA^\natural/\sN^\natural$ is pseudo-abelian,
where $\sN^\natural$ is
the ideal $\sN$ of $\sA^\natural$; but the obvious functor $\sA/\sN\to
\sA^\natural/\sN^\natural$ is clearly a pseudo-abelian envelope.
\end{proof}

\begin{thm}\label{t2} If $\bar\sA$ is of homological origin, any $\bar M\in \bar \sA$ is Schur-finite; in particular, $\bar\sA$ is
geometrically of integral type.
\end{thm}

\begin{comment}
\begin{proof} By Lemma \ref{l5}, we may assume $K$ algebraically closed.
Choose $(\sA,H:\sA\to
\Vec_L^\pm)$ such that $\bar A\simeq\sA/\sN$. Without loss of
generality, we may also assume that $L$ is algebraically closed. The
functor $H$ canonically
extends to a realisation
$H_L:\sA_L\to \Vec_L^\pm$. By \cite[Lemme 1]{err-nrsm},
$\sN(\sA_L)=\sN\otimes_K L$. Denote the
functor of extension of scalars $\bar\sA\to \bar\sA_L$ by $M\mapsto M_L$.
Clearly,
$\mu(M_L)=\mu(M)\otimes_K 1$, and the simple factors of $\End_{\bar
\sA_L}(M_L)=\End_{\bar \sA}(M)\otimes_K L$ are the same as those of
$\End_{\bar \sA}(M)$. All
this reduces us to the case where $K=L$ is algebraically closed. Without
loss of generality we
may further assume $H$ to be faithful. Finally, Lemma \ref{l7} reduces us
to the case
where
$\bar\sA$ is pseudo-abelian. We then have to prove that $\mu(S)\in\Z$ for
any simple $S\in \bar
\sA$.

Let $\tilde S$ be an object of $\sA$ which maps to $S$. By hypothesis,
$\End_{\bar
\sA}(S)\allowbreak=K$. Thus $\mu(S)\in K$, hence
$\tr(1_S)=\mu(S)$. On the other hand,
\begin{equation}\label{eq1}
\tr(1_S)=\tr(1_{\tilde S})= \dim_\gr H(\tilde S)\in \Z.
\end{equation}
\end{proof}
\end{comment}

\enlargethispage*{20pt}

\begin{proof} Choose $(\sA,H:\sA\to
\Vec_L^\pm)$ as in the proof of Lemma \ref{l5}. Without loss of
generality, we may assume that $H$ is faithful. Lift $\bar M$ to $M\in \sA$. Then $H(M)$ is
finite-dimensional, hence Schur-finite, which implies that $M$ and therefore $\bar M$ is
Schur-finite. The conclusion now follows from Proposition
\ref{p2} e).
\end{proof}

\begin{rk} The converse of Theorem \ref{t2} holds: namely, if every object of $\bar\sA$ is
Schur-finite, then the same is true in $(\bar\sA_{\bar K})^\natural$. By \cite[0.6 and
following remark]{deligne}, $(\bar\sA_{\bar K})^\natural$ is $\otimes$-equivalent to
$Rep(G,\epsilon)$, where $(G,\epsilon)$ is a super-affine group scheme over $\bar K$; in
particular, $(\bar\sA_{\bar K})^\natural$ and hence $\bar\sA$ admits a realisation into
$\Vec_{\bar K}$. This shows that \emph{if $\bar \sA$ is of homological origin, then it is
actually of homological type}. Another approach to this idea is the one in
\cite{aknote}, using $\otimes$-sections.

However, in the case of pure motives, one wants of course to study the general
situation of Definition \ref{d3} b), which is the one that arises naturally! This is what we do
in the remainder of this section.
\end{rk}

\begin{defn}\label{d2} Let $\sA$ be of homological type, and let
  $H:\sA\to\Vec_L^\pm$ be a realisation functor. Given $M\in \sA$, we
  say that \emph{the sign conjecture holds for $M$} (with respect to
  $H$) if there exists $p\in \End_\sA(M)$ such that $H(p)$ is the
  identity on $H^+(M)$ and is $0$ on $H^-(M)$.
\end{defn}

\begin{lemma}[cf. \protect{\cite[9.2.1]{nrsm}}] With the notation of
  Definition \ref{d2}:\\
a) If $M$ is finite-dimensional, it verifies the sign conjecture.\\
b) The converse is true if $H$ is faithful and $\sN(M,M)$ is a
nilideal (this is always the case if $L=K$).\qed
\end{lemma}

\begin{prop}\label{p5} Let $\sA$
be of homological type, $\bar\sA:=\sA/\sN$ and let $H:\sA\to
\Vec_L^\pm$ be a realisation 
functor. Then\\
a) For any simple object $S\in (\bar\sA_L)^\natural$, $d(S)\mid \mu(S)$.\\ 
b) Suppose $H$ faithful. Let $M\in
\sA$ verify the sign conjecture. Then the nilpotence level $r$ of
$\sN(M,M)$ verifies
\[r< \prod_i(\frac{|\mu_i(M)|}{e_i(M)}+1)\]
where $e_i(M)=d(S_i)$, with $S_i$ a simple summand of $\bar
M\in\bar\sA^\natural$ corresponding to its $i$-th isotypical
component (see Remark \ref{l1}).\\
c) If $\bar M$ is isotypical and $\mu(M)=\pm 1$, then $\sN(M,M)=0$. 
\end{prop}

\begin{proof} a) Since $\sA_L/\sN(\sA_L)=\bar\sA_L$, $\bar\sA_L$ is
neutrally of homological
origin; up to quotienting $\sA_L$ and replacing $K$ by $L$, we may assume
$L=K$ and $H$
faithful. Then $\sA$ is semi-primary and, as in the proof of Theorem
\ref{t2}, we may further
assume that $\sA$ and $\bar\sA$ are pseudo-abelian.

Let $\tilde S\in \sA$ mapping to $S$. By Wedderburn's theorem, the map
$\End_\sA(\tilde S)\to
\End_{\bar \sA}(S)$ has a ring-theoretic section $\sigma$. This makes
$H(\tilde S)$ a module
over the division ring $\End_{\bar \sA}(S)$. Therefore $\dim_K
H^\epsilon(\tilde S)$ is
divisible by $\dim_K \End_{\bar \sA}(S)=\delta(S)d(S)^2$ for $\epsilon =
\pm 1$. On the other
hand,
\[\dim_K H^+(\tilde S)-\dim_K H^-(\tilde S)=\mu(S)\delta(S)d(S)\]
by Proposition \ref{p2} a). Therefore, $\delta(S)d(S)^2$ divides
$\mu(S)\delta(S)d(S)$, which
means that $d(S)$ divides $\mu(S)$, as claimed.

b) Assume first $L=K$. Without loss of generality, we may also assume $\sA$
pseudo-abelian. Let 
$\sN=\sN(M,M)$ and consider the filtration $(\sN^i
H(M))_{0\le i\le
r-1}$. Note that $\sN^i H(M)=\sN^{i+1} H(M)$ $\iff$
$\sN^i=0$ since $\sN$ is a
nilpotent set of endomorphisms of $H(M)$. The associated graded
$(\gr^i H(M))_{0\le i\le r-1}$ is a graded $\End_{\bar
\sA}(\bar M)$-module, and $\gr^i
H(M)\ne 0$ for all $i<r$. 

Since $M$ verifies the sign conjecture, for each primary summand $\bar
M_i$ of $\bar M$, with lift $M_i$ in
$\sA$, one has either $H^+(M_i)=0$ or $H^-(M_i)=0$.  The proof of a) then
shows that $\gr H(M_i)$ is an $\End_{\bar \sA}(\bar
M_i)$-module of length $\frac{|\mu(S_i)|}{d(S_i)}$ where $S_i$ is the
associated simple object. Note that $\End_{\bar \sA}(\bar
M)=\prod_i\End_{\bar \sA}(\bar
M_i)$:  it follows that $\gr H(M)$ is an $\End_{\bar
\sA}(\bar M)$-module of length $\prod (\frac{|\mu(S_i)|}{d(S_i)}+1)-1$.
Hence the inequality.

In general we extend
  scalars from $K$ to $L$. Let
  $s=\prod_i(\frac{|\mu_i(M)|}{e_i(M)}+1)-1$. Applying 
  the result to the category $\sA L$ with the same objects as $\sA$
  and such that $\sA L(M,N)=H(\sA(M,N))L\subseteq \Hom_L(H(M),H(N))$,
  we get $(\sN(M,M)L)^s=0$, hence $\sN(M,M)^s\subseteq (\sN(M,M)L)^s=0$.

c) This follows immediately from b). 
\end{proof}

\begin{rk} In case $M$ is finite dimensional, we have another bound
  for the nilpotence level of $\sN(M,M)$ (valid without assuming $H$
  faithful). For simplicity, suppose that 
  $M$ is either positive or negative, and let $n = |\chi(M)|$. By
  \cite[Corollary 10.2]{kim}, $f^n=0$ for all $f\in \sN(M,M)$, which
  implies that $\sN(M,M)^{n^2}=0$ by a theorem of Razmyslov improving
  the Nagata-Higman bound \cite[4.3]{pro}. I learned of this
  improvement from Alessio del Padrone, who has generalised this bound
  to any finite-dimensional $M$ with $n=|\chi(M^+)|+|\chi(M^-)|$ \cite[2.4.10]{delp}.

  These two bounds have completely different behaviours: in Proposition
  \ref{p5} c) the former is optimal while the latter is not, but on
  the other hand if $M$ is a direct sum of $n$ invertible objects 
  pairwise non-isomorphic, the bound of Proposition
  \ref{p5} b) is $2^n-1$ while the other one is $n^2$.
  
As del Padrone pointed out, a third unrelated nilpotency bound is the one predicted by the
Bloch-Beilinson-Murre conjecture for the Chow motive of a smooth projective variety $X$
(namely, $\dim(X)+1$ \cite[Strong conj. 2.1]{jannsen2}).
\end{rk}

%\enlargethispage*{20pt}

\begin{rk} Coming back to the zeta function of an endomorphism, suppose
that $\sA$ is of
homological type; let $M\in \sA$ and $f\in \End_\sA(M)$. If $H$ is a
realisation of $\sA$, we
have by the usual computation
\[Z(f,t) = \det(1-H(f)t)^{-1}= \frac{\det(1-ft\mid H^-(M))}{\det(1-ft\mid
H^+(M))}.\]

Let $\bar M$ be the image of $M$ in $\bar \sA=\sA/\sN$ and $\bar f$ be
the image of $f$ in
$\End_{\bar \sA}(\bar M)$. Since $Z(f,t)=Z(\bar f,t)$, we get from
Theorem \ref{t2} and Theorem
\ref{t3} a) the identity
\[\frac{\det(1-ft\mid H^+(M))}{\det(1-ft\mid H^-(M))}=\prod_i
\Nrd_{A_i}(e_i-e_i \bar f
t)^{\mu_i(M)}.\]

Suppose for example that $\bar M$ is simple; the identity reduces to
\[\frac{\det(1-ft\mid H^+(M))}{\det(1-ft\mid H^-(M))}=\Nrd_{A}(1- \bar f
t)^{\mu(M)}\]
where $A=\End_{\bar \sA}(\bar M)$.

Supposing further that $\mu(M)>0$ to fix ideas, we find that the inverse
characteristic
polynomial of $f$ acting on $H^-(M)$ (with coefficients in $L$)
\emph{divides} the one for
$H^+(M)$, and the quotient has coefficients in $K$. This does not imply,
however, that
$H^-(M)=0$.
\end{rk}

\section{Examples}

In our first example, let $\sA$ be the rigid $K$-category of vector
bundles over $\P^1_K$.
It is of homological type, with realisation functor $H:\sA\to \Vec_L$
for $L=K(t)$ given by the generic fibre. Its indecomposable objects
are the $\sO(n)$ for $n\in\Z$: they are all of multiplicity $1$ but
$\sA(\sO(p),\sO(q))=\sN(\sO(p),\sO(q))\ne 0$ whenever $p<q$. This
shows that the 
condition ``isotypical'' is necessary in Proposition \ref{p5} c) (I am
indebted to Yves Andr\'e for pointing out this example). If one
extends scalars from $K$ to $L$ in the style of the proof of this
proposition, one finds $\sA L(\sO(p),\sO(q))=L$ whenever $p<q$.

Our main source of examples is, of course, the category $\sM_\num(k)$ of
pure motives over a field $k$ modulo numerical equivalence. By
Jannsen's theorem \cite{jannsen} and Theorem \ref{t2}, $\sM_\num(k)$
is semi-simple and geometrically of integral type. We shall compute
the multiplicities in certain cases.

We leave it to the reader to check that the multiplicities of Artin
(even Artin-Tate) motives are always $+1$. The next case is that of
abelian varieties. 

%\enlargethispage*{20pt}

Let $A$ be an abelian variety of dimension $g$ over $k$. Then we have the
Chow-K\"unneth decomposition
\[h(A)\simeq \bigoplus_{i=0}^{2g} h_i(A)\]
with $h_i(A)\simeq S^i(h_1(A))$ \cite{kunnemann}. Moreover,
\[\End h_1(A) =\End^0(A):= \End(A)\otimes\Q.\]

Moreover, $\chi(h_1(A))=-2g$. From this and Proposition \ref{p2}, we get
for $A$ simple:
\[\mu(A)=-\frac{2g}{\delta(\End^0(A))d(\End^0(A))}.\]

We recover the fact that the denominator divides the
numerator. 

Like Milne \cite{milne2}, we shall say that $A$ has \emph{many
  endomorphisms} if $\sum_i
\delta(\End^0(A_i))d(\End^0(A_i))=2g$, where $A_i$ runs through the
simple factors of $A$, or equivalently if
  $\delta(\End^0(A_i))d(\End^0(A_i))=2g_i$ for all $i$, where $g_i=\dim
  A_i$. This terminology is less ambiguous than ``having complex
  multiplication''. 

\begin{defn} Let $M\in \sM_\num(k)$ be a pure motive modulo numerical
  equivalence. Then $M$ is \emph{of abelian type} if it  is isomorphic to a direct summand of
the tensor product of an Artin-Tate motive and the motive of an abelian variety.
\end{defn}

Motives of abelian type are stable under direct sums, direct summands,
tensor products and duals. We then have:

\enlargethispage*{20pt}

\begin{thm}\label{tab} a) For $A$ a simple abelian variety,
  $\mu(h_1(A))=-1$ if and only if
  $A$ has many endomorphisms.\\
b) If $g=1$, then $\mu(h_1(A))=
\begin{cases}
-1&\text{if $A$ has complex multiplication}\\
-2&\text{otherwise}.
\end{cases}$\\
c) If $A$ has many endomorphisms, all multiplicities of $h_i(A)$
are equal to $(-1)^i$.\\
d) If $k$ is a finite field, then the multiplicities of any motive of
abelian type are equal to $\pm 1$.
\end{thm}

\begin{proof} a) and b) are clear; c) follows from a) and Corollary
  \ref{c2}, and d) follows from c) since any abelian variety over a
  finite field has many endomorphisms \cite{tate}.
\end{proof}

The next interesting case is that of $t_2(S)$ where $S$ is a surface
\cite{kmp}. If $k=\C$, there are many examples where the Hodge
realisation of $t_2(S)$ is simple \cite[Ex. 5 and Prop. 14]{ps}. A
fortiori $t_2(S)$ is simple, and Proposition \ref{p2} shows that  its
multiplicity equals its Euler characteristic, i.e. $b^2-\rho$ where $b^2$ is the
second Betti number and $\rho$ is the Picard number.

\section{An abstract version of the Tate(-Beilinson) conjecture}\label{5}

\subsection{Automorphisms of the identity functor }Let $\sA$ be a rigid
$K$-category, and let
$F$ be an \emph{$\otimes$-endomorphism of the identity functor of $\sA$}. By
\cite[I.5.2.2]{saa}, $F$ is then an isomorphism. Concretely, $F$ is given
by an automorphism
$F_M\in \End(M)$ for every object $M\in\sA$; $F_M$ is natural in
$M$, and further:
\begin{align*}
F_{M\oplus N} &= F_M\oplus F_N\\
F_{M\otimes N} &= F_M\otimes F_N\\
F_{M^*}&={}^t F_M^{-1} \text{ (\cf \cite[I, (3.2.3.6)]{saa}).}
\end{align*}

\begin{defn} The \emph{zeta function} (relative to $F$) of an object
$M\in\sA$ is
\[Z_F(M,t)=Z(F_M,t).\]
\end{defn}

\begin{lemma}\label{l8} The zeta function is additive in $M$:
\[Z_F(M\oplus N,t)=Z_F(M,t)Z_F(N,t).\]
It is multiplicative in $M$ in the following sense:
\[Z_F(M\otimes N,t)=Z_F(M,t)* Z_F(N,t)\]
where $*$ is the unique law on $1+tK[[t]]$ such that, identically,
$f*(gh)=(f*g)(f*h)$ and 
\[(1-at)^{-1} * (1-bt)^{-1} = (1-abt)^{-1}.\]
(Explicitly: if $f(t) = \exp\left(\sum_{n\ge 1} a_n\frac{t^n}{n}\right)$
and $g(t) = \exp\left(\sum_{n\ge 1} b_n\frac{t^n}{n}\right)$, then
$f*g(t) =  \exp\left(\sum_{n\ge 1} a_nb_n\frac{t^n}{n}\right)$.)\\
If moreover $\sA$ is semi-simple of integral type, then
\begin{enumerate}
\item $Z_F(M,t)\in K(t)$ for any $M\in \sA$;
\item $Z_F(M^*,t^{-1}) = (-t)^{\chi(M)}\det(F_M)Z_F(M,t)$;
\item for $S$ simple, 
\[Z_F(S,t) = P_S(t)^{-\chi(S)/\deg(F_S)}\]
where $P_S(t)$ is the inverse minimum polynomial of $F_S$ over $K$ and
$\deg(F_S)=\deg(P_S)=[K[F_S]:K]$.
\end{enumerate}
\end{lemma}

\begin{proof} Additivity is obvious; multiplicativity follows
from the identities
\[\tr(F_{M\otimes N}^n) =\tr(F_M^n\otimes F_N^n)=\tr(F_M^n)\tr(F_N^n).\]

(1), (2) and (3) follow from Theorem \ref{t3}: (1) from part a), (2) from
part b) by noting that $Z({}^tF_S^{-1},t^{-1})=Z(F_S^{-1},t^{-1})$, and
(3) from part a) again by noting that $F_S$ is in the centre of
$\End_\sA(S)$ (use Proposition \ref{p2} a)).
\end{proof}

\subsection{The semi-simple case}

\begin{defn} In the above, suppose $\sA$ semi-simple of integral type. We
say that $(\sA,F)$ \emph{verifies the Tate conjecture} if, for any $M\in
\sA$, $K[F_M]$ is the centre of $\End_\sA(M)$.
\end{defn}

\begin{thm}[\cf \protect{\cite[Th. 2.7]{geisser}}] \label{t4} Let $\sA$
be a semi-simple rigid pseudo-abelian $K$-category
of integral type, and let $F\in \Aut^\otimes(Id_\sA)$.  Then the
following conditions are equivalent: 
\begin{thlist}
\item  Given a simple object $S\in \sA$, $F_S=1_S$ implies $S=\1$.
\item  For any $M\in\sA$, $\ord_{t=1} Z_F(M,t) = -\dim_K \sA(\1,M)$. 
\item  For $S,T\in \sA$ simple, $P_S=P_T$ $\Rightarrow$ $S\simeq T$.
\item  For $M,N\in \sA$, $Z_F(M,t)=Z_F(N,t)$ $\Rightarrow$ $M\simeq N$.
\item $(\sA,F)$ verifies the Tate conjecture.
\end{thlist}
Moreover, these conditions imply:
\begin{itemize}
\item[(vi)] For any simple $S$, $|\mu(S)|=1$ and  $K[F_S]$ is the centre
of the algebra $\End_\sA(S)$.
%\item[(vii)] For $S,T\in \sA$ simple, $\End_\sA(S)\simeq \End_\sA(T)$
%$\Rightarrow$
%$S\simeq T$.
\end{itemize}
\end{thm}

\begin{proof} We shall prove the following implications:

(i) $\Rightarrow$ (ii)  $\Rightarrow$ (iii) $\Rightarrow$ (iv)
$\Rightarrow$ (i)

(ii)  $\Rightarrow$ (vi)

(iii) + (vi)  $\Rightarrow$ (v) $\Rightarrow$ (iii).

%(iii) + (v) $\Rightarrow$ (vii). 

(i) $\Rightarrow$  (ii): both sides are additive in $M$ so we may assume
$M$ simple. If $M=\1$, $Z_F(M,t)=1/(1-t)$ and the formula is true. If
$M\ne 1$, Lemma \ref{l8}
(3) and the hypothesis show that $\ord_{t=1}Z_F(M,t)=0$ and the formula
is also true.

 (ii) $\Rightarrow$  (iii):  Consider $f(t)=Z_T(S^*\otimes T,t)$. By Lemma
\ref{l8}, Formulas (2), (3) and the multiplicativity rule, we see that
\[f(t)=\prod_{i,j}(1-\alpha_i\alpha_j^{-1}t)^m\]
where $m=-\frac{\chi(S)}{\deg(F_S)}\frac{\chi(T)}{\deg(F_T)}$ and the
$\alpha_i$ are the roots
of the irreducible polynomial $P_S =P_T$ in a suitable extension of $K$.
Note that (ii) implies
that $\ord_{t=1} Z_F(M,t)\le 0$; the above formula shows that this
integer is $<0$. Hence
$0\ne\sA(\1,S^*\otimes T)\simeq \sA(S,T)$ and $S\simeq T$ by Schur's lemma.

(iii) $\Rightarrow$ (iv): write $M= \bigoplus_{i\in I} S_i^{m_i}$ and $N=
\bigoplus_{i\in I}
S_i^{n_i}$, where $S_i$ runs through a set of representatives of the
isomorphism classes of
simple objects of $\sA$. We then have, by Lemma \ref{l8} (3):
\begin{align*}
Z_F(M,t)& = \prod_{i\in I}P_{S_i}(t)^{-m_i\chi(S_i)/\deg(F_{S_i})}\\
Z_F(N,t) &= \prod_{i\in I}P_{S_i}(t)^{-n_i\chi(S_i)/\deg(F_{S_i})}.
\end{align*}

By hypothesis, the $P_{S_i}(t)$ are pairwise distinct irreducible
polynomials with constant
term $1$; then $Z_F(M,t)=Z_F(N,t)$ implies $m_i=n_i$ for all $i$, hence
$M\simeq N$.

(iv) $\Rightarrow$ (i): by hypothesis and Lemma \ref{l8} (3),
$Z_F(S,t)=(1-t)^{-\chi(S)}$. Thus
$Z_F(S,t)=Z_F(\1,t)^{\chi(S)}$. If $\chi(S)<0$, this gives
$S^{-\chi(S)}\simeq \1$, which
implies $\chi(S)=-1$ and $S\simeq \1$, which is absurd since
$\chi(\1)=1$. Thus $\chi(S)\ge 0$,
hence $S\simeq \1^{\chi(S)}$, hence $S\simeq \1$ since $S$ is simple.

(ii) $\Rightarrow$ (vi): the same computation as in the proof of (i)
$\Rightarrow$ (iii) gives
\begin{multline*}
\delta(S)d(S)^2=\dim\End_\sA(S) =-\ord_{t=1} Z(S^*\otimes S,t)\\
=\left(\frac{\chi(S)}{\deg(F_S)}\right)^2\ord_{t=1}
\prod_{i,j}(1-\alpha_i\alpha_j^{-1}t)=\frac{\chi(S)^2}{\deg(F_S)}.
\end{multline*}

Using the identity $\chi(S)=\mu(S)d(S)\delta(S)$ (\cf Proposition
\ref{p2} a)), we get
\[\deg(F_S)=\delta(S)\mu(S)^2.\]

But $\deg(F_S)\mid \delta(S)$, hence $\delta(S)=\deg(F_S)$ and $\mu(S)^2=1$.

(iii) + (vi) $\Rightarrow$ (v): Let $M=\bigoplus_i S_i^{m_i}$ with
$m_i>0$ and the $S_i$
simple and pairwise nonisomorphic. Then
\[\End_\sA(M)=\prod_i M_{m_i}(\End_\sA(S_i))\]
hence the centre of $\End_\sA(M)$ is the product of the centres of the
$\End_\sA(S_i)$. By
(vi), each of these centres is generated by $F_{S_i}$; by (iii), the
$P_{S_i}$ are pairwise
distinct irreducible polynomials, hence the minimum polynomial of $F_M$
must be divisible by
their product.

(v) $\Rightarrow$ (iii) (compare \cite{geisser}): if $P_S=P_T$ but
$S\not\simeq T$, then
$\End_\sA(S\oplus T) =\End_\sA(S)\times \End_\sA(T)$, with centre
containing $L\times L$ for
$L=K[F_S]=K[F_T]$. But $F_{S\oplus T}$ is killed by $P_S=P_T$, a contradiction.
\end{proof}

\begin{rk}\label{r1} Condition (vi) is really weaker than the others: take
$F=1$ in $\sA$ the category of linear representations of a finite abelian
group over $K$
algebraically closed.
\end{rk}

\enlargethispage*{20pt}

\begin{prop}\label{p3} Let $\sA$ be semi-simple of integral type and let $F\in
\Aut^\otimes(Id_\sA)$.\\ 
a) The Tate conjecture is true for $(\sA,F)$ if and only if it is true
for $(\sA^\natural, F)$,
where $\sA^\natural$ is the pseudo-abelian envelope of $\sA$ and $F$ is
extended to
$\sA^\natural$ naturally.\\ 
b) If $\sA$ is geometrically of integral type, the Tate conjecture
is invariant under extension of sca\-lars: if $L$ is an extension of $K$, then
$(\sA,F)$ verifies the Tate conjecture if and only if $(\sA_L,F)$
verifies the Tate conjecture.
\end{prop}

\begin{proof} a) ``If" is obvious. For ``only if", let $M=(N,e)\in
\sA^\natural$ where $N\in
\sA$ and $e$ is an idempotent of $N$. Write $M= \bigoplus_{i\in I}
S_i^{m_i}$ and $N=
\bigoplus_{i\in I} S_i^{n_i}$ as in the proof of Theorem \ref{t4}, (iii)
$\Rightarrow$ (iv). We
have
\[
\End(M)=\prod_i M_{m_i}(\End(S_i)),\quad
\End(N)=\prod_i M_{n_i}(\End(S_i)).
\]

Letting $Z_i$ denote the centre of $\End(S_i)$, we get
\[
Z(\End(M))=\prod_{m_i>0}Z_i,\quad Z(\End(N))=\prod_{n_i>0}Z_i.
\]

By hypothesis, $Z(\End(N))$ is generated by $F_N$ as a $K$-algebra; this
implies that $Z_i$ is
generated by $F_{S_i}$ for all $i$ and that the $P_{S_i}$ are pairwise
distinct. Hence $F_M$
generates  $Z(\End(M))$ as well.

b) This is obvious since the centre of a semi-simple algebra behaves well
under extension of
scalars.
\end{proof}

\begin{cor} If $(\sA,F)$ verifies the Tate conjecture, then the conditions of
Theorem \ref{t4} hold in $\sA$ even if $\sA$ is not pseudo-abelian.
\end{cor}

\begin{proof} This is obvious except for (ii) and (iv); but by
Proposition \ref{p3} a),
$(\sA^\natural,F)$ verifies the Tate conjecture; by Theorem \ref{t4},
$\sA^\natural$ also
verifies conditions (ii) and (iv), which a fortiori hold in its full
subcategory $\sA$.
\end{proof}

%\enlargethispage*{20pt}

\begin{prop}\label{p4} Suppose that
$(\sA,F)$ verifies the Tate conjecture. Let $S\in \sA$ be a simple object.\\ 
a) If $\chi(S)\ge 0$, then
$\Lambda^{\chi(S)+1}(S)=0$; if $\chi(S)< 0$, then  
$\bS^{-\chi(S)+1}(S)\allowbreak=0$.\\ 
b) $\sA$ is finite-dimensional in the sense of Kimura-O'Sullivan; more
precisely, there exists a
unique $\otimes$-$\Z/2$-grading of $\sA$ such that $S$ simple is positive
(\resp negative) if
and only if $\chi(S)>0$ (\resp $<0$). 
\end{prop}

\begin{proof} a) By Theorem \ref{t4} (iv), it suffices to see that
$Z_F(N,t)=1$ for
$N=\Lambda^{\chi(S)+1}(S)$ (\resp $N=\mathbf{S}^{-\chi(S)+1}(S)$): this follows
from the
computations of \cite[7.2.4]{nrsm}. b) is an immediate consequence (see also
\cite[9.2.1]{nrsm}).
\end{proof}

\subsection{The homological case}

Let $\sA$ be of homological type, provided with a
faithful realisation
functor $H:\sA\to \Vec_L^\pm$. Let $F\in \Aut^\otimes(Id_\sA)$, and let
us still denote by $F$
its image in $\Aut^\otimes(Id_{\bar \sA})$, where $\bar \sA=\sA/\sN$.
Note that $F$ acts on $H$ by functoriality. 

\begin{thm}\label{t5} Consider the following
conditions on an object $M\in \sA$:
\begin{thlist}
\item $\bar M\in \bar \sA$ verifies Condition {\rm (ii)} of Theorem \ref{t4}.
\item The map $\sA(\1,M)\otimes_K L\to H(M)^F$ is
surjective and the composition $H(M)^F\to H(M)\to H(M)_F$ is an isomorphism
(semi-simplicity at $1$).
\item The map $\sA(\1,M)\otimes_K L\to H(M)^F$ is
surjective and
$\sN(\1,M)=0$.
\item The sign conjecture holds for $M$.
\item $H^-(M)^F=0$.
\end{thlist}
Then 
\begin{enumerate}
\item {\rm (i) + (v)} $\iff$ {\rm (ii) + (iii)}. 
\item {\rm (i) + (iv)} $\Rightarrow$ {\rm (v)} $\Rightarrow$ {\rm (iv)}. 
\item {\rm (ii)} for $M$ and $M^*$ $\iff$ {\rm (iii)} for $M$ and $M^*$.
\end{enumerate}
\end{thm}

\begin{proof} These are classical arguments that only need to be put
straight in this abstract context.

Note that $H^-(\1)=0$, so that $\sA(\1,M)\otimes_K L\to H(M)^F$ actually
lands into $H^+(M)^F$; denote its image by $\sA(\1,M)L$. By definition of
$\sN$, the projection $\sA(\1,M)\otimes_K L\to \bar \sA(\1,M)\otimes_K L$
factors through $\sA(\1,M)L$. The 
diagram
\[\begin{CD}
\sA(\1,M)L&\hookrightarrow& H^+(M)^F\\
@V{\text{surj}}VV\\
\bar\sA(\1,\bar M)\otimes_K L
\end{CD}\]
gives the inequalities
\[\dim_L H^+(M)^F\ge \dim_L \sA(\1,M)L\ge \dim_K\bar \sA(\1,\bar M).\]

On the other hand,
\begin{multline*}
\ord_{t=1} Z_F(M,t) =\\
 \ord_{t=1} \det(1-F_Mt\mid H^{-}(M)) - \ord_{t=1}
\det(1-F_Mt\mid H^{+}(M))\\
=\dim_L H^{-}(M)^{F^\infty}-\dim_L H^{+}(M)^{F^\infty}
\end{multline*}
where $H^\pm(M)^{F^\infty}$ denotes the characteristic subspace of
$H^\pm(M)$ for the
eigenvalue $1$ under the action of $F$.

%\enlargethispage*{20pt}

(1) Suppose that $H^{-}(M)^F=0$. Then $H^{-}(M)^{F^{\infty}}=0$ and,
under (i), we have
\begin{multline*}
\dim_L H^+(M)^F\ge \dim_L \sA(\1,M)L\ge \dim_K\bar \sA(\1,\bar M)\\
= \dim_L H^+(M)^{F^\infty}\ge\dim_L H^+(M)^F
\end{multline*}
hence we have equality everywhere, and (ii) and (iii) are true.
Conversely, (ii) + (iii) gives isomorphisms $\bar\sA(\1,M)_L\iso
H(M)^F\iso H(M)^{F^\infty}$. In particular, $H^{-}(M)^F=0$ and we have
$\dim_K\bar \sA(\1,M)=\dim_L H^{+}(M)^{F^\infty}$, hence (i) and (v). Thus, (i) + (v) $\iff$
(ii) + (iii).

(2) Under (iv), we may write $M=M^+\oplus M^-$, with $H(M^+)$ purely even
and $H(M^-)$ purely odd. To prove that $H^{-}(M)^F=0$, we may therefore
consider separately the cases where $M$ is even and odd.

If $M$ is even, this is obvious. If $M$ is odd, we get, under (i):
\[H^+(M)^F= \sA(\1,M)=\bar \sA(\1,\bar M)=0\]
since $\sA(\1,M)\hookrightarrow H^+(M)^F$, and
\[-\dim H^-(M)^{F^\infty} = \dim \bar \sA(\1,\bar M)=0\]
which shows that (i) + (iv) $\Rightarrow$ (v). For (v) $\Rightarrow$
(iv), we reason as in \cite[Proof of Th. 2]{km}: there exists a
polynomial $\Pi\in K[t]$ such that $\Pi$ is divisible by $P^{-}$ and
$\Pi-1$ is divisible by $P^+$, where $P^\epsilon(t) =
 \det(t-F\mid H^\epsilon(M))$; then $\Pi(F)\in \End(M)$ is such that
$H(\Pi(M))$ is the identity on $H^+(M)$ and is $0$ on $H^-(M)$.

(3) The counit map $M\otimes M^*\to \1$ gives compatible pairings
\begin{align*}
\bar \sA(\1,\bar M)\times \sA(\1,\bar M^*)&\to K\\
\sA(\1,M)L\times \sA(\1,M^*)L&\to L\\
H(M)\times H(M^*)&\to L.
\end{align*}

The first and last are perfect pairings: for the first, check it on
simple objects thanks to
Schur's lemma\footnote{Or use the definition of the ideal $\sN$.} and for
the last, this
follows from the structure of the tensor category $\Vec_L^\pm$. Consider
now the commutative
diagram
\[
\begin{CD}
\bar \sA(\1,\bar M)_L@<a<\text{surj}< \sA(\1,M)L@>b>\text{inj}> H(M)^F\\
@V{\wr}VV @VcVV @VdVV\\
(\bar \sA(\1,\bar M^*)_L)^*@>a^*>\text{inj}>(\sA(\1,M^*)L)^*@<b^*<\text{surj}<
(H(M^*)^F)^*&\simeq H(M)_F.
\end{CD}\]

Notice that the right vertical map coincides with the one of (ii). 

Now assume that $b$ and $b^*$ are isomorphisms. The diagram shows
immediately that $a,a^*$
isomorphisms $\Rightarrow$ $d$ isomorphism. Conversely, if $d$ is an
isomorphism, so is $c$;
but then, $a$ and $a^*$ must be isomorphisms. Finally, $a$ is an
isomorphism $\Rightarrow$
$\sA(\1,M)\to \bar \sA(\1,M)\otimes_K L$ is injective $\Rightarrow$
$\sN(\1,M)=0$, as desired. 
\end{proof}

%\enlargethispage*{20pt}

\begin{cor}[\cf \protect{\cite[2.9]{tate2}}]\label{c1} Let $\sA,H,F$ be as in Theorem \ref{t5}, and
suppose that $\sA$ is
pseudo-abelian. Consider the following conditions:
\begin{thlist}
\item The Tate conjecture holds for $(\bar\sA,F)$.
\item $\sA\to\bar \sA$ is an equivalence of categories and $H$ induces a
fully faithful functor
\[\tilde H:\bar \sA_L\to \Rep_L(F)_{ss}^\pm\]
where the right hand side denotes the $\otimes$-category of $\Z/2$-graded
$L$-vector spaces
provided with the action of an automorphism $F$, this action being semi-simple.
\item The sign conjecture holds for any $M\in \sA$ (equivalently
\cite[9.2.1 c)]{nrsm}, $\sA$
is a Kimura-O'Sullivan category).
\item For any $M\in \sA$, $H^{-}(M)^F=0$.
\end{thlist}
Then {\rm (iv)} $\Rightarrow$ {\rm (iii)}; moreover {\rm (ii)} $\iff$
{\rm (i) + (iii)} $\iff$
{\rm (i) + (iv)}.\\ 
If these conditions are verified, then for any simple object $S\in
\sA_L^\natural$, $\End(S)$ is commutative.
\end{cor}

\begin{proof} First, (iv) $\Rightarrow$ (iii) by Point 2 of Theorem \ref{t5}.
 If now (ii) holds, then Conditions (ii) and (iii) of Theorem \ref{t5}
hold for any
$M$, hence so do its conditions (i) and (v) by Point 1 of this theorem.
Point 2 also shows
that $M$ verifies Condition (iv) of this theorem. This shows that (ii)
$\Rightarrow$ (i) +
(iii) + (iv) in Corollary \ref{c1}.

Suppose that (i) holds. Then Condition (ii) of Theorem \ref{t4} holds for
any $\bar M\in
\bar\sA$. If moreover $H^{-}(M)^F=0$ for any $M\in \sA$, Conditions (ii)
and (iii) of Theorem
\ref{t5} are verified for any $M\in \sA$ by Point 1 of this theorem.
Applying this to
$M=P^*\otimes Q$ for some $P,Q\in \sA$, the adjunction isomorphisms
\[\sA(P,Q)\simeq \sA(\1,P^*\otimes Q)\]
show that $\sN(P,Q)=0$, hence a bijection
\[\bar\sA(M,N)\otimes_K L\to \Hom_F(H(M),H(N)).\]

Moreover, since $\bar\sA$ is semi-simple, $H(F_M)$ is a semi-simple
endomorphism of $H(M)$ for
any $M\in \sA$. This shows that (i) + (iv) $\Rightarrow$ (ii).   

Suppose that (i) and (iii) hold. Then Point 2 of Theorem \ref{t5} shows that
$H^{-}(M)^F=0$ for any $M\in \sA$, thus (i) + (iii) $\Rightarrow$ (iv).

It remains to justify the last claim: it follows from Proposition
\ref{p5} and Condition (vi)
of Theorem \ref{t4}.
\end{proof}

\begin{rk} In the classical case of motives over a finite field,
Conditions (iii) and (iv) of Corollary \ref{c1} hold
provided the Weil cohomology $H$ verifies the Weak Lefschetz theorem, by
Katz-Messing
\cite{km}. It is a little annoying not to be able to dispense of them in
this abstract setting,
especially in view of Proposition \ref{p2} b).
\end{rk}

\subsection{The Tate-Beilinson conjecture}

We conclude by transposing the argument of \cite{cell} to this
abstract context.

\begin{thm}[\cf \protect{\cite[Th. 1]{cell}}]\label{tcell} Let $\sA$ be a rigid $K$-category
  provided with a 
  $\otimes$-auto\-morph\-ism $F$ of the identity functor. Suppose 
  that $\sN$ is locally nilpotent (e.g. that $\sA$ is a
  Kimura-O'Sullivan category), and that $\bar\sA=\sA/\sN$ verifies the
  Tate conjecture relatively to $F$. Then $\sN=0$, i.e. $\sA=\bar\sA$.
\end{thm}

\begin{proof} We note that the hypothesis on $\sN$ implies that the
  functor $\sA\to\bar\sA$ is conservative. The argument is the same as
  in \cite{cell}: by rigidity
  it is sufficient to show that $\sA(\1,M)\iso \bar\sA(\1,\bar M)$ for
  any  $M\in \sA$. By the nilpotence of $\sN(M,M)$, we may lift to
  $\End_\sA(M)$ an orthogonal system of idempotents of $\bar M$
  corresponding to a decomposition in simple summands. This reduces us
  to the case where $\bar M$ is simple. There are two cases:

\begin{enumerate}
\item $\bar M\simeq \1$. Then $M\simeq 1$ by conservativity, and both
  Hom groups are isomorphic to $K$.
\item $\bar M\not\simeq\1$. Then $\bar\sA(\1,\bar M)=0$ and we have to
  show that $\sA(\1,M)\allowbreak=0$. By Theorem \ref{t4} (i),
  $F_{\bar M}\ne 1$,
  hence by conservativity, $1-F_M$ is an isomorphism. If now $f\in
  \sA(\1,M)$, we have $f=F_M f$, hence $f=0$.
\end{enumerate}
\end{proof}

\enlargethispage*{20pt}

\end{document}